\documentclass[12pt]{article}

\usepackage{amssymb,latexsym,amsmath}

\oddsidemargin=\evensidemargin \footskip=36pt \voffset=-1.5cm
\newenvironment{proof}[1][Proof]{\textbf{#1.} }{\ \rule{0.5em}{0.5em}}
\def\thrm{\begin{theorem}}
\def\thrml#1{\begin{theorem}\label{#1}}
\def\ethrm{\end{theorem}}
\def\prpstn{\begin{proposition}}
\def\prpstnl#1{\begin{proposition}\label{#1}}
\def\eprpstn{\end{proposition}}
\def\rmrk{\begin{remark}}
\def\rmrkl#1{\begin{remark}\label{#1}}
\def\ermrk{\end{remark}}
\def\dfntn{\begin{definition}}
\def\dfntnl#1{\begin{definition}\label{#1}}
\def\edfntn{\end{definition}}
\def\nmrt{\begin{enumerate}}
\def\enmrt{\end{enumerate}}
\def\tm#1{\item[{\rm (#1)}]}
\def\qtn{\begin{equation}}
\def\qtnl#1{\begin{equation}\label{#1}}
\def\eqtn{\end{equation}}
\def\lmm{\begin{lemma}}
\def\lmml#1{\begin{lemma}\label{#1}}
\def\elmm{\end{lemma}}
\def\crllr{\begin{corollary}}
\def\crllrl#1{\begin{corollary}\label{#1}}
\def\ecrllr{\end{corollary}}
\def\hpthss{\begin{hypothesis}}
\def\hpthssl#1{\begin{hypothesis}\label{#1}}
\def\ehpthss{\end{hypothesis}}
\def\xmpl{\begin{example}}
\def\xmpll#1{\begin{example}\label{#1}}
\def\exmpl{\end{example}}
\def\css{\begin{cases}}
\def\ecss{\end{cases}}

\DeclareMathOperator{\aut}{Aut}
\DeclareMathOperator{\AGL}{AGL}
\DeclareMathOperator{\alt}{Alt}
\DeclareMathOperator{\AGaL}{A{\mathrm\Gamma}L}

\DeclareMathOperator{\cyc}{Cyc}

\DeclareMathOperator{\GL}{GL}

\DeclareMathOperator{\GF}{GF}
\DeclareMathOperator{\GaL}{\mathrm\Gamma L}
\DeclareMathOperator{\Inv}{Inv}
\DeclareMathOperator{\id}{id}

\DeclareMathOperator{\Iso}{Iso}
\DeclareMathOperator{\Mat}{Mat}
\DeclareMathOperator{\orb}{Orb}

\DeclareMathOperator{\rk}{rk}
\DeclareMathOperator{\SL}{SL}
\DeclareMathOperator{\SP}{Sp}
\DeclareMathOperator{\SU}{SU}
\DeclareMathOperator{\soc}{soc}
\DeclareMathOperator{\sym}{Sym}

\def\F{{\mathbb F}}
\def\K{{\mathbb K}}

\def\N{{\mathbb N}}
\def\Z{{\mathbb Z}}

\def\CC{{\cal C}}

\def\NN{{\cal N}}

\def\R{{\cal R}}
\def\S{{\cal S}}

\def\zv{0_{\scriptscriptstyle{V}}}

\def\ov{\overline}

\def\proof{{\bf Proof}.\ }
\def\bull{\vrule height .9ex width .8ex depth -.1ex }

\newtheorem{formula}{}[section]
\newtheorem{proposition}[formula]{Proposition}
\newtheorem{definition}[formula]{Definition}
\newtheorem{corollary}[formula]{Corollary}
\newtheorem{remark}[formula]{Remark}
\newtheorem{lemma}[formula]{Lemma}
\newtheorem{theorem}[formula]{Theorem}
\newtheorem{problem}[formula]{Problem}
\newtheorem{hypothesis}[formula]{Conjecture}
\newtheorem{example}[formula]{Example}

\def\thrm{\begin{theorem}}
\def\thrml#1{\begin{theorem}\label{#1}}
\def\ethrm{\end{theorem}}
\def\prpstn{\begin{proposition}}
\def\prpstnl#1{\begin{proposition}\label{#1}}
\def\eprpstn{\end{proposition}}
\def\rmrk{\begin{remark}}
\def\rmrkl#1{\begin{remark}\label{#1}}
\def\ermrk{\end{remark}}
\def\dfntn{\begin{definition}}
\def\dfntnl#1{\begin{definition}\label{#1}}
\def\edfntn{\end{definition}}
\def\nmrt{\begin{enumerate}}
\def\enmrt{\end{enumerate}}
\def\tm#1{\item[{\rm (#1)}]}
\def\qtn{\begin{equation}}
\def\qtnl#1{\begin{equation}\label{#1}}
\def\eqtn{\end{equation}}
\def\lmm{\begin{lemma}}
\def\lmml#1{\begin{lemma}\label{#1}}
\def\elmm{\end{lemma}}
\def\crllr{\begin{corollary}}
\def\crllrl#1{\begin{corollary}\label{#1}}
\def\ecrllr{\end{corollary}}
\def\hpthss{\begin{hypothesis}}
\def\hpthssl#1{\begin{hypothesis}\label{#1}}
\def\ehpthss{\end{hypothesis}}
\def\xmpl{\begin{example}}
\def\xmpll#1{\begin{example}\label{#1}}
\def\exmpl{\end{example}}
\def\prblm{\begin{problem}}
\def\prblml#1{\begin{problem}\label{#1}}
\def\eprblm{\end{problem}}
\def\css{\begin{cases}}
\def\ecss{\end{cases}}

\makeatletter
\renewcommand{\subsection}{\@startsection{subsection}{2}{0mm}{-2mm}{-2mm}
{\bf\normalsize}}
\def\sbsn{\subsection{\hspace{-5mm}}}

\makeatother

\title{ On cyclotomic schemes over finite near-fields}
\author{
J. Bagherian \\
Institute of Advanced Studies in Basic Sciences\\[-3pt]
P.O.Box:\ 45195-1159, Zanjan, Iran \and
Ilia Ponomarenko
\thanks{Partially supported by RFFI 05-01-00899, NSH-4329.2006.1.
The author thanks the Institute for Advanced Studies in Basic Sciences
(IASBS) for its hospitality during the time that this paper was
written.}\\[-1pt]
Petersburg Department of V.A.Steklov\\[-3pt]
Institute of Mathematics\\[-3pt]
Fontanka 27, St. Petersburg 191023, Russia\\[-3pt]
{\tt inp@pdmi.ras.ru}\\[-3pt]
http://www.pdmi.ras.ru/\~{}inp \and
A. Rahnamai Barghi
\thanks{Corresponding author: rahnama@iasbs.ac.ir } \\
Institute of Advanced Studies in Basic Sciences\\[-3pt]
P.O.Box:\ 45195-1159, Zanjan, Iran}

\begin{document}

\date{}
\maketitle

\begin{abstract}
We introduce a concept of cyclotomic association scheme $\CC$ over a finite near-field. It is
proved that if $\CC$ is nontrivial, then $\aut(\CC)\le\AGL(V)$ where $V$ is the linear space associated
with the near-field. In many cases we are able to get more specific information about $\aut(\CC)$.
\vspace{3pt}

\vspace{2mm}

\end{abstract}

\pagebreak\vspace{0in}

\section{Introduction}
An algebraic structure $\K=(\K,+,\circ)$ is called a (right) {\it near-field} if $\K^+=(\K,+)$
is a group with the neutral element $0_\K$, $\K^\times=(\K\setminus\{0_\K\},\circ)$ is a group,
$x\circ 0_\K=0_\K$ for all $x\in\K$ and
\qtnl{f280105a}
(x+y)\circ z=x\circ z+y\circ z,\qquad x,y,z\in \K.
\eqtn
In finite case the group $\K^+$ is elementary abelian and the group $\K^\times$ is abelian iff
$\K$ is a field (concerning the theory of near-fields we refer to \cite{Wa87}).
Moreover, by the Zassenhaus theorem apart from seven exceptional cases each finite near-field
$\K$ is the {\it Dickson near-field}, i.e. there exists a finite field $\F_0$ and its extension
$\F$ such that $\F^+=\K^+$ and
\qtnl{f080306b}
y\circ x=y^{\sigma_x}\cdot x,\qquad x,y\in\K
\eqtn
where $\sigma_x\in\aut(\F/\F_0)$ and $\cdot$ denotes the multiplication in $\F$.
In this case $|\F_0|=q$ and $|\K|=|\F|=q^n$ where $q$ is a power of a certain prime $p$
and $n=[\F:\F_0]$. It can be proved that $(q,n)$ forms a {\it Dickson pair} which means that every
prime factor of $n$ is a divisor of $q-1$ and $4\,|\, n$ implies $4\,|\, (q-1)$. There exist
exactly $\varphi(n)/k$ nonisomorphic Dickson near-fields corresponding to the same Dickson
pair $(q,n)$ where $k$ is the order of $p\,(\mathrm{mod}\, n)$.

Let $\K$ be a finite near-field and $K$ be a subgroup of the group $\K^\times$. Set
$\R=\{R_a\}_{a\in\K}$ where
\qtnl{f250405a}
R_a=\{(x,y)\in\K^2:\ y-x\in a\circ K\}
\eqtn
is a binary relation on the set $\K$. Then it is easily seen that any such a relation is a
2-orbit of the permutation group
\qtnl{f290106c}
\Gamma(K,\K)=\{x\mapsto x\circ b+c,\ x\in\K:\ b\in K,c\in\K\}
\eqtn
and so the pair $(\K,\R)$ forms an {\it association scheme} on $\K$ (see Section~\ref{f280405a} for
the background on permutation groups and association schemes). We call it the {\it cyclotomic
scheme} over the near-field $\K$ and denote it by $\cyc(K,\K)$. The number $|K|$ is called
the {\it valency} of the scheme. If $K=\K^\times$, then the scheme is of rank $2$ and we
call it the {\it trivial} scheme. The set of all cyclotomic schemes of valency $m<q^n-1$ over a
Dickson near-field corresponding to a Dickson pair $(q,n)$, is denoted by $\cyc(q,n,m)$.

When $\K=\F$ is a field, we come to cyclotomic schemes introduced by P.~Delsarte (1973). One
can see that any two such schemes of the same rank are isomorphic. Moreover, the automorphism
group of such a nontrivial scheme is a subgroup of the group $\AGaL_1(\F)$ (see
\cite[p.389]{BCN}). However there exist a number of cyclotomic schemes over near-fields which are
not isomorphic to cyclotomic schemes over fields. The main purpose of this paper is to study
isomorphisms of cyclotomic schemes over near-fields.

The additive group of a finite near-field $\K$ being an elementary abelian one,
can be identified with the additive group of a linear space $V_\K$ over the prime field
containing in the center of $\K$.
If there exist isomorphic cyclotomic schemes over near-fields $\K$ and $\K'$, then
obviously $|\K|=|\K'|$ and hence the linear spaces $V_\K$ and $V_{\K'}$ are isomorphic. Thus
to study isomorphisms of cyclotomic schemes we can restrict ourselves to near-fields $\K$
with fixed linear space $V=V_\K$.

\thrml{f250705a}
Let $\CC$ and $\CC'$ be nontrivial cyclotomic schemes over near-fields $\K$ and $\K'$
respectively. Suppose that $V=V_\K=V_{\K'}$. Then $\Iso(\CC,\CC')\subset\AGL(V)$. In particular,
$\aut(\CC)\le\AGL(V)$.
\ethrm

For the trivial scheme $\CC$ we obviously have
$\aut(\CC)=\sym(\K)$. Thus the inclusion $\aut(\CC)\le\AGL(V)$ holds only if $|\K|\le 4$.
In general case, the right-hand side of the first inclusion of Theorem~\ref{f250705a}
can not be reduced because $\Iso(\CC,\CC)=\AGL(V)$
where $\CC=\cyc(K,\F)$ with $\F$ being a finite field of composite order and $K$ being the
multiplicative group of the prime subfield of $\F$.

We prove Theorem~\ref{f250705a} in Section~\ref{f200305g}. The key point of the proof is
Theorem~\ref{f110405a} showing that the operation of taking the 2-closure preserves the socle of
any uniprimitive 3/2-transitive permutation groups of the affine type. (Here we essentially
use the result of \cite{PS}.) This also gives a criterion for the isomorphism of nontrivial
cyclotomic schemes (Theorem~\ref{f260405a}).
\vspace{2mm}

The second part of Theorem~\ref{f250705a} can be made more precise in some cases. For instance,
if the cyclotomic scheme $\CC=\cyc(K,\K)$ is imprimitive, then $\aut(\CC)=\Gamma(K,\K)$
(Corollary~\ref{f210305d}). In general case it is not true even for a cyclotomic scheme over
a finite field because $\aut(\CC)$ can contain some automorphisms of the field. However,
we are able to restrict the automorphisms of a cyclotomic scheme by using Zsigmondy prime
divisors of its valency.

\dfntn
Given integers $q,n\in\N$ a prime divisor $r$ of $q^n-1$ is called a {\it Zsigmondy prime} for
$(q,n)$ if $r$ does not divide $q^i-1$ for all $1\le i<n$. The set of all such primes
greater than a natural number $k$ is denoted by $Z_k(q,n)$.
\edfntn

It is known that at least one Zsigmondy prime for $(q,n)$ exists unless $(q,n)=(2,6)$ or $q+1$ is a power of
2 and $n=2$ (see e.g.~\cite{R97}). Moreover, any such prime is of the form $r=an+1$ for some
$a\ge 1$.

\thrml{f301205a}
Let $\CC\in\cyc(p^d,n,m)$ be a cyclotomic scheme over a Dickson near-field. Suppose that $m$ has
a prime divisor $r\in Z_{2dn+1}(p,dn)$. Then the group $\aut(\CC)$ is isomorphic to a subgroup
of the group $\AGaL_1(p^{dn})$.
\ethrm

From Corollary~\ref{f290106a} it follows that in the condition of Theorem~\ref{f301205a}
the scheme $\CC$ is primitive. In fact, this theorem shows also that for a sufficiently large $n$
the group $\aut(\CC)$, is isomorphic to a subgroup of the group $\AGaL_1(p^d)$ in all but one case.

\thrml{f301205b}
Let $\CC\in\cyc(p^d,n,m)$ be a cyclotomic scheme over a Dickson near-field. Then for $n\gg q=p^d$
the group $\aut(\CC)$ is isomorphic to a subgroup of the group $\AGaL_1(q^n)$ unless
$Z_{2dn+1}(p,dn)=\{r\}$ with $r^2\nmid (q^n-1)$ where $r=(q^n-1)/m$.
\ethrm

Let $\CC=\cyc(K,\K)$ and $\CC'=\cyc(K',\K)$ be cyclotomic schemes of the same valency
$m<q^n-1$ over a Dickson near-field $\K$ of order $q^n$. Then
$$
[\K^\times :K]=\rk(\CC)-1=(q^n-1)/m=\rk(\CC')-1=[\K^\times :K'].
$$
Suppose that $r=(q^n-1)/m$ is a prime and $r^2\nmid (q^n-1)$. Then the groups $K$ and $K'$
are the Hall subgroups of the group $\K^\times$ and so are conjugate in it. So the schemes
$\CC$ and $\CC'$, and hence the groups $\aut(\CC)$ and $\aut(\CC')$, are isomorphic.
Thus Theorem~\ref{f301205b} shows that given a Dickson near-field $\K$ corresponding to
the Dickson pair $(q,n)$ with
$n\gg q$, there is at most one (up to isomorphism)
nontrivial cyclotomic scheme $\CC$ over $\K$ for which we don't know whether $\aut(\CC)$ is
isomorphic to a subgroup of the group $\AGaL_1(q^n)$.

Theorem~\ref{f301205b} is proved in Section~\ref{f221105a} by means of the classification of linear
groups with orders having certain large prime divisors, given in~\cite{GPPS}. We believe that
more delicate analysis of this classification could improve our result to show that
$\aut(\CC)$ is isomorphic to a subgroup of the group $\AGaL_1(q^n)$ apart from a finite
number of possible $n$'s for a fixed~$q$.

\section{Permutation groups and association schemes}\label{f280405a}

\sbsn Concerning basic facts of finite permutation group theory we refer to~\cite{DM}. For a
positive integer $m$ and a group $\Gamma\le\sym(V)$ the set of orbits of the induced action of
$\Gamma$ on $V^m$ is denoted by $\orb_m(\Gamma)$; the elements of this set are called {\it $m$-orbits}
of $\Gamma$. The group $\Gamma$ is called {\it $m$-closed} iff it coincides with its {\it $m$-closure}
$\Gamma^{(m)}$ which is by definition the largest subgroup of $\sym(V)$ with the same set of
$m$-orbits as $\Gamma$.

Let $U$ be a set with at least two elements and $m\ge 2$ be an integer. Following \cite{PS} we
say that a permutation group $G\le\sym(V)$ {\it preserves a product decomposition} $U^m$ of
$V$, if the latter can be identified with the Cartesian product $U^m$ in such a way that $G$ is
a subgroup of the wreath product $\sym(U)\wr\sym(m)$ in product action. Any element $g$ of the
latter group induces uniquely determined permutations $g_1,\ldots,g_m\in\sym(U)$ and
$\sigma=\sigma_g\in\sym(m)$ such that
\qtnl{f270305c}
(u_1,\ldots,u_m)^g=(u_{i_1}^{g_{i_1}},\ldots,u_{i_m}^{g_{i_m}})\quad
\mbox{where}\quad i_j=j^{\sigma^{-1}}.
\eqtn
If $G$ projects onto a transitive subgroup of $\sym(m)$, then the subgroup of index~$m$
in $G$ stabilizing the first entry of points of $U^m$ induces a subgroup of $\sym(U)$ on the set
$U$ of first entries of points of $V=U^m$; this group is called the group {\it induced by $G$ on
$U$}. The following statement being a special case of result of~\cite[Lemma~4.1]{PS} will be
used in Section~\ref{f200305g}. Below a primitive group is called {\it uniprimitive} if
it is not 2-transitive, and it is called of {\it affine type} if the socle of it is abelian.

\thrml{f080405a}
Let $G\le\sym(V)$ be a uniprimitive group of the affine type. Suppose that
$\soc(G)\ne\soc(G^{(2)})$. Then $G$ and $G^{(2)}$ preserve a product decomposition $V=U^m$ such
that $|U|\ge 5$, $m\ge 2$ and the group induced by $G^{(2)}$ on $U$ contains $\alt(U)$.\bull
\ethrm

\sbsn
Let $V$ be a finite set and $\R$ be a partition of the set $V^2$ containing the diagonal
$\Delta(V)$ of $V^2$. Then the pair $\CC=(V,\R)$ is called a (noncommutative) {\it association
scheme} or a {\it scheme} on $V$ if $\R$ is closed with respect to the permutation of coordinates,
and given the binary relations $R,S,T\in\R$ the number
$$
|\{v\in V: (u,v)\in R, (v,w)\in S\}|
$$
does not depend on the choice of $(u,w)\in T$. Two schemes $\CC=(V,\R)$ and $\CC'=(V',\R')$ are
called {\it isomorphic} if $\R^f = \R'$, for some bijection $f:V\to V',$ called the
{\it isomorphism} from $\CC$ to $\CC'$, where $\R^f =\{ R^f\in\R: R\in\R\}$ and
$R^f=\{(u^f,v^f):(u,v)\in R\}.$ The set of all such isomorphisms is denoted by $\Iso(\CC,\CC')$.
The group $\Iso(\CC)=\Iso(\CC,\CC)$ contains a normal subgroup
$$
\aut(\CC)=\{g\in\sym(V):\ R^g = R,\ R\in \R\}
$$
called the {\it automorphism group} of $\CC$.

The elements of $V$ and $\R=\R(\CC)$ are called the {\it points} and the {\it basis relations}
of~$\CC$ respectively; the numbers $\deg(\CC)=|V|$ and $\rk(\CC)=|\R|$ are called the {\it degree}
and the {\it rank} of~$\CC$. The scheme $\CC$ is called {\it imprimitive} if there exists an
equivalence relation $E$ on $V$ such that $E\not\in\{\Delta(V),V^2\}$ and $E$ is a union of some
basis relations of $\CC$; otherwise $\CC$ is called {\it primitive} whenever $\deg(\CC)>1$.

A wide class of schemes comes from permutation groups as follows. Let $\Gamma\leq\sym(V)$ be a
permutation group and $\R$ the set of all 2-orbits of $\Gamma$. Then the pair
$\Inv(\Gamma)=(V,\R)$ is a scheme and
$$
\aut(\Inv(\Gamma))=\Gamma^{(2)}.
$$
In particular, any cyclotomic scheme $\cyc(K,\K)$ over a near-field $\K$ equals the scheme
$\Inv(\Gamma)$ with $\Gamma=\Gamma(K,\K)$. It is also true that the scheme $\Inv(\Gamma)$
is primitive iff so is the group $\Gamma$.

\sbsn
Let $\K$ be a near-field and $K\le\K^\times$. Then the group $\Gamma(K,\K)$ defined by
(\ref{f290106c}) can be naturally identified with a subgroup of the group $\AGL(V)$ where
$V=V_\K$. Under this identification the group $K$ (considered as the subgroup of the group
$\Gamma(K,\K)$) goes to a subgroup of the group $\GL(V)$. This subgroup is called
the {\it base group} of the corresponding cyclotomic scheme $\cyc(K,\K)$. Thus the base
group is nothing else than the image of the natural linear representation of $K$ in $\GL(V)$.

\thrml{f280705a}
Let $\CC$ be a cyclotomic scheme over a near-field $\K$. Then $\CC$ is primitive iff the base
group of $\CC$ is irreducible.
\ethrm
\proof Let $\CC=\cyc(K,\K)$ for some group $K\le\K^\times$. Then $\CC=\Inv(\Gamma)$ where
$\Gamma=\Gamma(K,\K)$. So the scheme $\CC$ is primitive iff the group $\Gamma$ is primitive.
However, from \cite[p.19]{Su} it follows that the latter statement holds iff the stabilizer of
the point $0_\K$ in the group $\Gamma$ is an irreducible subgroup of the group $\GL(V_\K)$. Since
this stabilizer coincides with the base group of the scheme $\CC$, we are done.\bull
\vspace{2mm}

It should be noted that for a primitive cyclotomic scheme $\cyc(K,\K)$ the base
group can be primitive (as a linear group) or not. The latter case is realized e.g. for $\K$
being the field of order~$9$ and for $K$ being the subgroup of $\K^\times$ of order~$4$ (then
$\Gamma(K,\K)$ is isomorphic to the subgroup of the wreath product $\sym(3)\wr\sym(2)$ in the
product action).

\crllrl{f290106a}
Let $\CC$ be a cyclotomic scheme satisfying the hypothesis of Theorem~\ref{f301205a}. Then
$\CC$ is primitive.
\ecrllr
\proof Let $\CC=\cyc(K,\K^\times)$ and $G$ be the base group of $\CC$. Then $G\le\GL(V)$ where
$V=V_\K$ and $r$ divides $|G|=|K|=m$. So the group $G$ is irreducible (see~\cite[\S~6]{MW}) and
we are done by Theorem~\ref{f280705a}.\bull
\vspace{2mm}

Let $V$ be a linear space over a prime finite field and $G\le\GL(V)$ be an irreducible abelian
group. Then $G$ is cyclic and the linear span $\F_0=L(G)$ of it (in the algebra $\Mat(V)$) is a finite
field with $|V|$ elements (see~\cite{MW}). In particular, the group $\F_0^\times$ being a Singer
subgroup of $\GL(V)$, acts regularly on the set $V^{\#}=V\setminus\{0_V\}$ by the multiplication
of matrix to vector. So for a fixed element $u\in V^\#$ there is a bijection
$$
\tau:\F_0\to V,\quad f\mapsto u^f
$$
translating the field structure from $\F_0$ to $V$. For the corresponding field $\F=\F(G)$ we
have $\F^+=V^+$ and $K\le\F^\times$ where $K$ is the permutation group on $V$ induced by
the group~$G$.

\thrml{f020505a}
Any primitive cyclotomic scheme with the abelian base group is a cyclotomic scheme over a field.
\ethrm
\proof Let $\CC=\cyc(K,\K)$ be a primitive cyclotomic scheme over a near-field $\K$ where the
group $K\le\K^\times$ is abelian. The base group $G\le\GL(V)$ where $V=V_\K$, of this
scheme is isomorphic to $K$ and hence is also abelian. Due to the primitivity of $\CC$ from
Theorem~\ref{f280705a} it follows that $G$ is irreducible. By the definition of the field
$\F=\F(G)$ we have $\F^+=V^+=\K^+$ and $K\le\F^\times\cap\K^\times$. Thus
$\CC=\Inv(\Gamma(K,\K))=\Inv(\Gamma(K,\F))$ is a cyclotomic scheme over the field $\F$ and we
are done.\bull

\section{An isomorphism criterion for cyclotomic schemes}\label{f200305g}

\sbsn
In this section we prove Theorem~\ref{f250705a}. When the base group of a cyclotomic scheme is
primitive as a linear group, the required statement immediately follows from
Theorem~\ref{f080405a}. In the imprimitive case we need to strengthen the latter theorem by
means of the following lemma. We recall that a transitive group $\Gamma\le\sym(V)$ is called
{\it 3/2-transitive} if all the orbits of its one point stabilizer $\Gamma_v$ on $V\setminus\{v\}$
have the same size.

\lmml{f060405a}
Let $G\le\sym(V)$ be a 3/2-transitive group preserving a product decomposition $V=U^m$ where
$m\geq 2$. Then $G_{u,v}$ is an abelian 2-group for distinct points $u,v\in V$.
\elmm
\proof Let us fix a point $u=(u_0,\ldots,u_0)\in U^m$ where $u_0\in U$. Then from (\ref{f270305c})
it follows that $u_0^{g_i}=u_0$ for all $g\in G_u$ and all $i\in I=\{1,\ldots,m\}$. This implies
that the cardinality of the set $I_v=\{i\in I:\ v_i\ne u_0\}$ with $v_i$ being the $i$th component
of $v$, does not depend on the choice of $v$ inside of an orbit of the group $G_u$. Thus, the sets
$$
V_k=\{v\in V:\ |I_v|=k\},\qquad k\in\Z,
$$
\qtnl{f080405b}
R=\{(v,w)\in V_1\times V_2:\ v_i=w_i\quad\text{for all}\quad i\in I_v\cap I_w\}
\eqtn
are $G_u$-invariant. We note that from the definition of $R$ it follows that $|R_{in}(w)|= 2$
for all $w\in V_2$ where $R_{in}(v)=\{u\in V:\ (u,v)\in R\}$.
\vspace{2mm}

{\bf Claim 1.} Let $X_1\in\orb(G_u,V_1)$, $X_2\in\orb(G_u,V_2)$ and $S=R_{X_1,X_2}$ (here and
below we set $R_{X,Y}=R\cap(X\times Y)$ for all $X,Y\subset V$). Then
$$
|S_{out}(x)|= 2,\qquad x\in X_1,
$$
where $S_{out}(u)=\{v\in V:\ (u,v)\in S\}$. Indeed, since $S$ is a $G_u$-invariant relation the
numbers $|S_{out}(x)|$ and $|S_{in}(v)|$ do not depend on $x\in X_1$ and $v\in X_2$
respectively. Denote them by $a_1$ and $a_2$. Then $|X_1|a_1=|X_2|a_2$. Taking into account that
$|X_1|=|X_2|$ due to 3/2-transitivity of $G$, we conclude that $a_1=a_2$. Since $a_2=2$ by the
definition of the relation $S$ (see (\ref{f080405b})), and we are done.
\vspace{2mm}

{\bf Claim 2.} The following inequality holds:
$$
|y^{G_{u,x}}|= 2,\qquad x,y\in V_1,\quad I_x\ne I_y.
$$
Indeed, let $x,y\in V_1$. Then $I_x=\{i\}$ and $I_y=\{j\}$ for some distinct $i,j\in I$. So
there exists the uniquely determined element $v\in V_2$ such that $x_i=v_i$ and $y_j=v_j$. Then
$(x,v),(y,v)\in R$. From Claim~1 with $X_1$ and $X_2$ being the orbits of the group $G_u$
containing $x$ and $v$, it follows that $S_{out}(x)=\{v,v'\}$ for some $v'\in X_2\backslash\{v'\}$
where $S=R_{X_1,X_2}$. It is easy to see that the set $S_{out}(x)$ is $G_{u,x}$-invariant and
hence so is the set $R_{in}(v)\cup R_{in}(v')$. However, this set contains at most three points
and two of them are $x$ and $y$. So
$$
|y^{G_{u,x}}|\le |(R_{in}(v)\cup R_{in}(v'))\setminus\{x\}|= 2
$$
and we are done.
\vspace{2mm}

{\bf Claim 3.} Let $x\in V_1$ and $v\in V_2$. Then $(G_{u,x})_Y$ is a 2-group where
$Y=v^{G_{u,x}}$. Indeed, let $I_v=\{i,j\}$ for some $i,j\in I$. Since $i\ne j$, we can assume
that $\{i\}\ne I_x$. Set $y$ to be the unique element of $V_1$ such that $y_i=v_i$. Then from
Claim~2 it follows that $y^{G_{u,x}}=\{y,z\}$ for some $z\in V_1$. So
$$
Y\subset S_{out}(y)\cup S_{out}(z)
$$
where $S=R_{X,Y}$ with $X=y^{G_{u,x}}$. Moreover, by Claim~1 we also have that both of sets in
the right-hand side are of cardinality equal~2. Thus taking into account that
$S_{in}(v)\cap\{y,z\}=\{y\}$ we see that either $|Y|=4$ and $|(G_{u,x,y,z})_Y|=2$. In any case
$|(G_{u,x})_Y|\in\{4\}$ and we are done.
\vspace{2mm}

{\bf Claim 4.} The action of $G_u$ on $V_2$ is faithful. Indeed, let $g\in G_u$ be such that
$v^g=v$ for all $v\in V_2$. Take $v\in V_2$ with $I_v=\{i,j\}$ and $v_i=v_j=u'$ where $u'\in
U\setminus\{u_0\}$. Then it follows that $(u')^{g_i}=(u')^{g_j}=u'$. This implies that
$g_i=\id_U$ for all $i\in I$. On the other hand, taking $v_i\ne v_j$, we see that
$i^{\sigma_g}=i$ and $j^{\sigma_g}=j$. Thus $\sigma=\id_I$ and we are done.
\bull
\vspace{2mm}

To complete the proof of Lemma~\ref{f060405a} take $v\in V_1$. By Claim 4 the group $G_{u,v}$
acts faithfully on $V_2$. So it is isomorphic to a subgroup of the direct product of the groups
$(G_{u,v})_X$ where $X$ runs over the set $\orb(G_{u,v},V_2)$. Then it is 2-group by
Claim~3 and we are done.\bull

\thrml{f110405a}
Let $G\le\sym(V)$ be a uniprimitive 3/2-transitive group of the affine type. Then
$\soc(G)=\soc(G^{(2)})$.
\ethrm
\proof Suppose that $\soc(G)\ne\soc(G^{(2)})$. Then from Theorem~\ref{f080405a} it follows that
the groups $G$ and $G^{(2)}$ preserve a product decomposition $V=U^m$ such that $|U|\ge 5$,
$m\ge 2$ and the group induced by $G^{(2)}$ on $U$ contains $\alt(U)$. This implies that
\qtnl{f110405b}
|G^{(2)}|=am|\alt(U)|
\eqtn
for some natural number $a$. On the other hand, the group $G^{(2)}$ is obviously uniprimitive
and 3/2-transitive. So the size of any nontrivial orbit of a one point stabilizer of $G^{(2)}$
equals to the same number, say $d$. One can see that $d=me$ for some divisor $e$ of
$|U|-1$. So by Lemma~\ref{f060405a} applied to $G^{(2)}$ we have
\qtnl{f110405c}
|G^{(2)}|=|V|me2^k
\eqtn
for some natural number $k$. Now from (\ref{f110405b}) and (\ref{f110405c}) it follows that
$$
|U|\frac{(|U|-1)}{e}\frac{(|U|-2)!}{2}\quad\mbox{divides}\quad |V|2^k.
$$
However, this is impossible for $|U|\ge 5$. Indeed, $|V|=p^b=|U|^m$ for some prime $p$ and
some natural number $b$, but for
$|U|\ge 5$ the number in the left-hand side has at least one prime divisor different
from $p$ and~2.\bull
\vspace{2mm}

Since the 2-closure of 3/2-transitive group is 3/2-transitive, from Theorem~\ref{f110405a}
it follows that the group $G^{(2)}$ is a uniprimitive 3/2-transitive group of the affine type.
If in addition, $G$ preserves a product decomposition, then the same decomposition is
preserved by $G^{(2)}$. Thus the form of this group can be found by means of the classification
of 3/2-transitive imprimitive linear groups given in~\cite{Ps68}.

\sbsn
In this subsection we fix a near-field $\K$, a cyclotomic scheme $\CC$ over $\K$ and
denote by $T=T_V$ the translation group of the linear space~$V=V_\K$. In particular,
$T\le\sym(V)$.

\lmml{f210305a}
If $\rk(\CC)>2$, then $T$ is a characteristic subgroup of the group~$\aut(\CC)$. More
exactly,
\nmrt
\tm{1} if $\CC$ is imprimitive, then $\aut(\CC)$ is a Frobenius group with the Frobenius kernel
$T$,
\tm{2} if $\CC$ is primitive, then $T=\soc(\aut(\CC))$.
\enmrt
\elmm
\proof Set $\Gamma=TG$ where $G$ is the base group of $\CC$. Since $\aut(\CC)=\Gamma^{(2)}$,
the orbits of the stabilizer of the point $\zv$ in the group $\aut(\CC)$ coincide with the
orbits of the group $G$. On the other hand, obviously,
$$
|X|=|G|,\qquad X\in\orb(G,V^\#)
$$
where $V^\#=V\setminus\{\zv\}$. So $\aut(\CC)$ is a 3/2-transitive permutation group. Suppose
first that the scheme $\CC$ is imprimitive. Then the group $\aut(\CC)$ is imprimitive. Since
any 3/2-transitive group is either primitive or a Frobenius group \cite[Theorem 10.4]{W64},
we conclude that $\aut(\CC)$ is a Frobenius group. The Frobenius kernel of this
group has the cardinality $|V|=|T|$ and contains all fixed point free elements of the group
$\Gamma$. Thus the Frobenius kernel coincides with $T$ which proves statement~(1).

Suppose that $\CC$ is a primitive scheme. Then the group $\Gamma$ is primitive. Clearly, the
group $T$ is a normal abelian subgroup of it. This implies that the socle of $\Gamma$ is
abelian and hence $\soc(\Gamma)=T$ (see \cite[Theorem~4.3.B]{DM}). Thus $\Gamma$ is a group of
the affine type. Since it is 3/2-transitive, Theorem~\ref{f110405a} implies that
$$
T=\soc(\Gamma)=\soc(\Gamma^{(2)})=\soc(\aut(\CC))
$$
which completes the proof.\bull
\vspace{2mm}

{\bf Proof of Theorem~\ref{f250705a}.} Let $f\in\Iso(\CC,\CC')$. Then obviously $f$ is a
permutation group isomorphism from $\aut(\CC)$ to $\aut(\CC')$. Since both of these
groups are transitive, without loss of generality we assume that $f$ leaves the point $0_V$
fixed. Then it suffices to verify that $f\in\aut(T)=\GL(V)$. However, the schemes $\CC$ and
$\CC'$ are primitive or not simultaneously. Thus the required statement follows from
Lemma~\ref{f210305a}.\bull

\sbsn
To make statements of Theorem~\ref{f250705a} more precise given a group $G\le\GL(V)$ we set
\qtnl{f020505b}
\ov G=G^{(1)}\cap\GL(V).
\eqtn
Clearly, $\ov G$ coincides with the largest group $H\le\GL(V)$ such that $\orb(H)=\orb(G)$.
We call this group the {\it linear closure} of $G$.

\thrml{f260405a}
In the conditions of Theorem~\ref{f250705a} denote by $G$ and $G'$ the base groups of the
schemes $\CC$ and $\CC'$ respectively. Then these schemes are isomorphic iff the groups $\ov G$
and $\ov {G'}$ are conjugate in $\GL(V)$. Moreover, $\aut(\CC)=T\ov G$ where $T=T_V$.
\ethrm
\proof From Lemma~\ref{f210305a} it follows that $T$ is a normal subgroup of the group
$\aut(\CC)=(TG)^{(2)}$. Thus the second statement of the theorem is the consequence of
(\ref{f020505b}) and the following lemma.

\lmml{f270405a}
Let $A$ be a group and $G\le\aut(A)$. Denote by $T$ the permutation group on $A$ induced
by the right regular representation of $A$. Suppose that $T$ is a normal subgroup of the
group $\Gamma=(TG)^{(2)}$. Then $\Gamma=T(\aut(A)\cap G^{(1)})$.
\elmm
\proof Set $H$ to be the stabilizer of the point $e=1_A$ in the group $\Gamma$. Then
$\orb(H)=\orb(G)$ and hence $H\le G^{(1)}$. On the other hand, since $T$ is normalized by $H$
and $e^H=\{e\}$ we have $h^{-1}t_ah=t_{a^h}$ for all $a\in A$ and $h\in H$ where $t_a$ is
the element of $T$ taking $x$ to $xa$. So
$$
(ab)^h=a^{t_bh}=(a^h)^{h^{-1}t_bh}=(a^h)^{t_{b^h}}=a^hb^h,\qquad a,b\in A.
$$
Thus $H\le\aut(A)$, and hence $\Gamma=TH\le T(\aut(A)\cap G^{(1)})$. Conversely, let
$g\in\aut(A)\cap G^{(1)}$. Then given $c\in A$ there exists $g_c\in G$ such that $c^g=c^{g_c}$.
So given $a,b\in A$ due to the normality of $T$ we have
$$
(a,b)^g=(a^t,b^t)^{t^{-1}g}=(e,c)^{gs}=(e^g,c^g)^s=(e,c^{g_c})^s=(e^{g_c},c^{g_c})^s=
(a^{tg_c},b^{tg_c})^s=(a,b)^{tg_cs}
$$
where $t$ is the element of $T$ such that $a^t=e$, $c=b^t$ and $s=g^{-1}t^{-1}g$.
Since $tg_cs\in TG$, this means that $g$ preserves the 2-orbit of the group $TG$ containing
$(a,b)$ for all $a,b\in A$. Thus $g\in\Gamma$ and hence, $\aut(A)\cap G^{(1)}\le\Gamma$. This
implies that $T(\aut(A)\cap G^{(1)})\le\Gamma$ and we are done.\bull
\vspace{2mm}

To prove the first statement of Theorem~\ref{f260405a} suppose that $g^{-1}\ov G g=\ov{G'}$ for
some $g\in\GL(V)$. Then by the second statement of the theorem we have
$$
g^{-1}\aut(\CC)g=g^{-1}(T\ov G)g=(g^{-1}Tg)(g^{-1}\ov Gg)=T\ov{G'}=\aut(\CC')
$$
whence it follows that $g\in\Iso(\CC,\CC')$. Conversely, let
$g\in\Iso(\CC,\CC')$. Then $g^{-1}\aut(\CC)g=\aut(\CC')$. By
Theorem~\ref{f250705a} without loss of generality we can assume
that $g\in\GL(V)$. Then $g$ leaves the point $v=\zv$ fixed and by
the first part we have
$$
g^{-1}\ov Gg=g^{-1}(\aut(\CC))_vg=\aut(\CC')_v=\ov{G'}
$$
which complete the proof.\bull
\vspace{2mm}

For imprimitive cyclotomic schemes Theorem~\ref{f260405a} can be simplified as follows.

\crllrl{f210305d}
Let the cyclotomic schemes $\CC$ and $\CC'$ be imprimitive. Then they are isomorphic
iff their base groups are conjugate in $\GL(V)$. Moreover, $\ov G=G$ and $\aut(\CC)=TG$.
\ecrllr
\proof From statement (1) of Lemma~\ref{f210305a} it follows that $\aut(\CC)$ is an imprimitive
Frobenius group. Since $TG\leq\aut(\CC)$ and $|TG|= |\aut(\CC)|$, it follows that $\aut(\CC)=TG$.
Since obviously $G\le\ov G$, this also shows that $\ov G=G$. Now the first part of the
required statement is a consequence of Theorem~\ref{f260405a} after taking into account that
two isomorphic schemes are primitive or not simultaneously.\bull

\section{Proof of Theorems~\ref{f301205a} and \ref{f301205b}}\label{f221105a}

The main tool of this section is the following theorem which is deduced from
the classification \cite{GPPS} of linear groups with orders having certain large prime divisors.
In our case such a divisor coincides with a Zsigmondy prime $r$ for a pair $(q,n)$. We
observe that any cyclic group $G\le\GL(n,q)$ of order $r$ is irreducible. This is a
consequence of the fact that the linear span $L(G)$ of it is a finite field $\F$ with $q^n$
elements. Below we consider the group $\GaL_1(\F)$ as a subgroup of $\GL(n,q)$.

\thrml{f290106f}
Let $G\le\Gamma\le\GL(n,q)$ where $(q,n)\not\in\{(2,4),(2,6)\}$. Suppose that $G$ is
a cyclic group of order $r\in Z_{2n+1}(q,n)$ and the group $\Gamma$ acts intransitively on the
set of all nonzero vectors of the underlying linear space. Then $\Gamma\le\GaL_1(\F)$ where
$\F=L(G)$.
\ethrm
\proof We observe that
the Zsigmondy prime $r$ for $(q,n)$ is a primitive prime divisor of $q^n-1$ in terms of
\cite{GPPS}. Since $r$ divides $|\Gamma|$, the hypothesis of the theorem implies that
the group $\Gamma$
satisfies the condition of the Main Theorem of that paper with $d=e=n$. So by this theorem
one of the following statements holds:
\nmrt
\tm{1} $\Gamma$ has a normal subgroup isomorphic to $\SL_n(q)$, $\SP_n(q)$, $\SU_n(q)$ or
$\Omega_n^-(q)$,
\tm{2} $r\le 2n+1$,
\tm{3} $\Gamma\le\GL(n/m,q^m)\cdot m$ for some divisor $m\ne 1$ of $n$,
\tm{4} $(q,n)=(2,4)$ or $(2,6)$
\enmrt
where $\GL(n/m,q^m)\cdot m$ is the general linear group $\GL(n/m,q^m)$ embedded to $\GL(n,q)$
and extended by the group of automorphisms of the field extension $\GF(q^m):\GF(q)$.
However, the cases (1), (2) and (4) are contradict to the intransitivity of $\Gamma$,
the conditions on $r$ and on $(q,n)$ respectively. Let us consider case (3).

Suppose first that $m\ne n$. Set $q'=q^m$ and $n'=n/m$. Then obviously
$(q',n')\not\in\{(2,4),(2,6)\}$ and $r\ge 2n+1>2n'+1$. This implies that
$Z_{2n+1}(q,n)\subset Z_{2n'+1}(q',n')$ and we can apply the same arguments by induction.
Thus we can assume that $m=n$ and $\Gamma\le\GL(n/n,q^n)\cdot n=\GaL_1(\F')$ for some field
$\F'\subset\Mat(n,q)$. To complete the proof we observe that the multiplicative group of $\F'$
is contained in the normalizer $\GaL_1(\F)$ of the group $G$ in $\GL(n,q)$. Since
$\F^\times$ is the unique Singer subgroup of $\GL(n,q)$ contained in $\GaL_1(\F)$ (see~\cite{CR}),
it follows that $\F'=\F$ and we are done.\bull
\vspace{2mm}

{\bf Proof of Theorem~\ref{f301205a}.} From the hypothesis of the theorem it follows that
$r$ divides the order $m$ of the base group of the scheme $\CC$. So the latter group contains
a cyclic subgroup $G$ of order $r$. By Theorem~\ref{f250705a} we see that $\Gamma=\aut(\CC)_u$
with $u=0_V$ is a subgroup of $\GL(nd,p)$. Moreover, $G\le\Gamma$ and $\Gamma$ acts
intransitively on the set $V^\#$ (because $m<p^d-1$). Finally, $(p,nd)\in\{(2,4),(2,6)\}$
only if $n=1$, because $(2,4),(4,2),(8,2)$ and $(2,6)$ are not Dickson pairs, and for
the Dickson pair $(4,3)$ we have $(p,d,n)=(2,2,3)$ and the set $Z_{2dn+1}(p,dn)$ is empty
(there are no Zsigmondy primes for $(2,6)$). Thus $\Gamma\le\GaL_1(\F)$ with $\F=L(G)$ by
Theorem~\ref{f290106f},
and hence the group $\aut(\CC)$ is isomorphic to a subgroup of the group $\AGaL_1(\F)$.
Since $|\F|=p^d$, we are done.\bull
\vspace{2mm}

{\bf Proof of Theorem~\ref{f301205b}.}
First we cite some number theoretical results from \cite{R00}. Given $n\in\N$ denote by $P[n]$
the greatest prime factor of $n$ and given $0<\kappa<1/\log 2$ let
$\NN_\kappa=\{n\in\N: D(n)\le\kappa\log\log n\}$ where $D(n)$ is the number of distinct prime
factors of $n$. Then according to \cite[p.25]{R00} given
real numbers $\alpha,\beta$ there exists a constant $C_\kappa>0$ such that for each $n\ge 3$,
$n\in\NN_\kappa$, the following inequality holds:
\qtnl{f290106g}
P[\Phi_n(\alpha,\beta)]>C_\kappa\frac{n\log n^{1-\kappa\log 2}}{\log\log\log n}
\eqtn
with $\Phi_n(\alpha,\beta)=\prod_i(\alpha-\zeta^i\beta)$ where $\zeta$ is a primitive $n$th
root of 1 and $i$ runs over the set of all numbers $1,\ldots,n$ coprime to~$n$.

Let us fix a prime power $q=p^d$ and given a number $N\in\N$ set
$$
\NN(q,N)=\{dn\in\N:\ dn\ge N\ \text{and}\ (q,n)\ \text{is a Dickson pair}\}.
$$
Choose $N_q$ to be the minimal number $n\in\N$ for which $D(d)\log q\le\kappa\log\log n$
where $\kappa=1/(2\log 2)$. Then given $dn\in\NN(q,N_q)$ we have
$$
D(dn)\le D(d)D(n)\le D(d)\log q\le \kappa\log\log n.
$$
So $dn\in\NN_\kappa$ and hence $\NN(q,N_q)\subset\NN_\kappa$. By (\ref{f290106g}) this implies
that for a fixed $q$ we have
\qtnl{f290106h}
P[\Phi_{dn}(p,1)]>C_\kappa\frac{n\sqrt{\log n}}{\log\log\log n}>2dn+1
\eqtn
for all sufficiently large $n\in\NN(q,N_q)$. On the other hand, a prime factor $r$ of the number
$\Phi_{dn}(p,1)$ is not a Zsigmondy prime for $(p,dn)$ iff $r\le dn$ (see
\cite[Proposition 2]{R97}). Thus from (\ref{f290106h}) it follows that there exists a natural
number $N'_q>N_q$ such that $Z_{2dn+1}(p,dn)\ne\emptyset$ for all $n>N'_q$.

To complete the proof let $\CC=\cyc(K,\K^\times)$ where $\K$ is a Dickson near-field
corresponding to the Dickson pair $(q,n)$. Suppose that $m=|K|<q^n$ and $n>N'_q$. Then
$Z_{2dn+1}(p,dn)\ne\emptyset$. Let us show that if
\qtnl{f290106i}
|Z_{2dn+1}(p,dn)|>1,\quad\text{or}\quad Z_{2dn+1}(p,dn)=\{r\}\quad\text{and}\quad r^2\,|\,(q^n-1)
\eqtn
where $r=(q^n-1)/m$, then the group $\aut(\CC)$ is isomorphic to a subgroup of the group
$\AGaL_1(q^n)$. To do this set $K'$ to be a maximal subgroup of $\K^\times$ containing $K$.
The group $\K^\times$ being isomorphic to a subgroup of the group $\GaL_1(q^n)$, is solvable.
Due to the maximality of $K'$ this implies that there exists a normal elementary abelian
subgroup $K_0$ of $\K^\times$ such that $[\K^\times:K']=|K_0|$ is a prime power. Moreover,
any Sylow subgroup of $\K^\times$ is a cyclic group or a quaternion group \cite{Wa87}.
Thus $K_0$ is a cyclic group of prime order and so the number $[\K^\times:K']$ is prime.
By (\ref{f290106i}) this implies that $m'=|K'|$ has a prime divisor $r'\in Z_{2dn+1}(p,dn)$.
So from Theorem~\ref{f301205a} applied to the cyclotomic scheme $\CC'=\cyc(K',\K^\times)$ it
follows that the group $\aut(\CC')$ is isomorphic to a subgroup of the group $\AGaL_1(q^n)$.
Since $\aut(\CC)\le\aut(\CC')$, we are done.\bull

\end{document}